\documentclass[a4paper,12pt]{smfarts}
\usepackage[left=1cm, right=1cm,paperheight=5in,paperwidth=7in]{geometry} 
\title
{Stable first-order theories as simplicial profinite sets} 
\thanks{Misha Gavrilovich.
Department of Mathematics, Ben-Gurion University of the Negev P.O.B. 653, Be’er Sheva
84105, Israel. {\tt{mi\!\!\!ishap\!\!\!p@sd\!\!\!df.org}}. \\
 This research was partially supported by the Israel   
 Science Foundation grant No. 2196/20. 
}
\usepackage{cfr-lm}
\usepackage{MnSymbol}
\usepackage{ae}

\usepackage{smfthm}
\usepackage{vmargin}
\usepackage{enumitem}
\usepackage{caption}
\usepackage{pdfpages}
\usepackage[matrix, arrow,all,cmtip,color]{xy}
\usepackage{url}
\usepackage{epigraph}

\def\includegraphics[#1]#2{    \pdfximage width  \linewidth {#2} \pdfrefximage\pdflastximage}
\newbox\TestBox
\def\Remove #1 {\setbox\TestBox=\hbox{#1}%
        \leavevmode\rlap{\vrule height 2.5pt depth-1.75pt width\wd\TestBox}%
	        \box\TestBox\ }

\def\includegraphicss[#1]#2{ 
	\pdfximage width  #1\linewidth {#2} 
	\pdfrefximage\pdflastximage 
}

\def\Imm{\operatorname{Im}}

\makeatletter
\newcommand{\colim@}[2]{%
	  \vtop{\m@th\ialign{##\cr
	      \hfil$#1\operator@font colim$\hfil\cr
	          \noalign{\nointerlineskip\kern1.5\ex@}#2\cr
		      \noalign{\nointerlineskip\kern-\ex@}\cr}}%
		      }
		      \newcommand{\colim}{%
			        \mathop{\mathpalette\colim@{\rightarrowfill@\scriptscriptstyle}}\nmlimits@
				}
	
\newcommand{\mylim@}[2]{%
	  \vtop{\m@th\ialign{##\cr
	      \hfil$#1\operator@font lim$\hfil\cr
	          \noalign{\nointerlineskip\kern1.5\ex@}#2\cr
		      \noalign{\nointerlineskip\kern-\ex@}\cr}}%
		      }
\newcommand{\mylim}{%
			        \mathop{\mathpalette\mylim@{\leftarrowfill@\scriptscriptstyle}}\nmlimits@
				}
\makeatother 	

\def\inv{^{-1}}

\def\Topp{\operatorname{Top}}

   \def\rrt#1#2#3#4#5#6{\xymatrix{ {#1} \ar[r]|{} \ar@{->}[d]|{#2} & {#4} \ar[d]|{#5} \\ {#3}  \ar[r] \ar@{-->}[ur]^{}& {#6} }}

\newcommand{\bi}{\begin{itemize}}
\newcommand{\ei}{\end{itemize}}

\def\bee{\begin{enumerate}[label={(\arabic*)},ref={(\arabic*)}]} 
\newcommand{\eee}{\end{enumerate}}

\def\bii{\begin{itemize}[label={(\arabic*)},ref={(\arabic*)}]} 
\newcommand{\eii}{\end{itemize}}

\def\biii{\begin{itemize}}
\newcommand{\eiii}{\end{itemize}}

\newcommand{\bd}{\begin{itemize}\item}
\newcommand{\ed}{\end{itemize}}

\def\calA{\mathcal A}

\def\BBBb{\Bbb B_\bullet}

\def\lra{\longrightarrow}
\def\rtt{\,\rightthreetimes\,}
\def\xra{\xrightarrow}

\usepackage{hyperref}
\usepackage{endnotes}

\def\rrt#1#2#3#4#5#6#7{\xymatrix{ {#1} \ar[r]^{} \ar@{->}[d]_{#2} & {#4} \ar[d]^{#5} \\ {#3}  \ar[r] \ar@{-->}[ur]^{#7}& {#6} }}

\def\rtt{\rightthreetimes}

\def\lra{\longrightarrow}
\def\lr{{\rtt lr}}
\def\lrl{{\rtt l}}

\def\rlr{{\rtt r}}

\def\id{\operatorname{id}}

\def\xra{\xrightarrow}

\def\Dop{\Delta^{\operatorname{op}}}

\def\PPhi{\ensuremath{\Filt}}

\newtheorem{que}[paragraph]{Question}

\usepackage[ethiop,english]{babel}
\usepackage{ethiop}

 \newcommand{\ethi}{\selectlanguage{ethiop}}
\def\Filt{{\ethi\ethmath{wA}}}

\def\FF{\ethi\ethmath{qa}}
\def\FFae{\ethi\ethmath{qo}}

\def\Funct{\text{Funct}}

\def\PPhi{\ensuremath{\Filt}}
\def\skiip#1{}

\def\smallblackbox{{\Huge\centerdot}}

\def\mathcalDense{\protect{\hphantom{\bullet_c\searrow}{\!\smallblackbox \atop  \displaystyle\Downarrow }\atop { 
\raisebox{3pt}{$\bullet$}\searrow\raisebox{-3pt}{$\smallblackbox$}  }}} 

\def\nondenseimage{\mathcalDense}

\def\cl{\operatorname{cl}}

\def\card{\operatorname{card}}
\def\pr{\operatorname{pr}}

\def\bqqq{\begin{quote}}
\def\eqqq{\end{quote}}

\newcommand{\forkindep}[1][]{%
  \mathrel{
    \mathop{
      \vcenter{
        \hbox{\oalign{\noalign{\kern-.3ex}\hfil$\vert$\hfil\cr
              \noalign{\kern-.7ex}
              $\smile$\cr\noalign{\kern-.3ex}}}
      }
    }\displaylimits_{#1}
  }
}
\usepackage{tikz}
\usepackage{tikz-cd}

\def\mathcalAnobrackets{{\displaystyle \bullet\,\,\bullet \atop{\displaystyle \displaystyle\Downarrow \atop\displaystyle\bullet}}}

\def\aand{\operatorname{\&}}

\def\BbbQ{\mathbb Q}
\def\SymQ{\operatorname{Sym}\BbbQ} 
\def\faktor#1#2{\left.\raisebox{1pt}{\ensuremath{#1}}\middle/\raisebox{-1pt}{\ensuremath{#2}}\right.}
\def\faktour#1#2{\faktor{#1\!}{\!#2}}

\def\FF{\ethi\ethmath{wA}}

\def\FFae{\ethi\ethmath{wI}}

\def\sFF{s\FF}

\def\sFFae{s\FFae}

\def\sPPhi{s\PPhi}

\def\BbbR{\Bbb R}

\def\Sing{\mathrm{Sing}}
\def\Singb{\Sing_\bullet}

\def\constb{\mathrm{const}_\bullet}

\def\calL{\mathcal L}

\def\Homm#1#2{{\mathrm{Hom}(#1,#2)}}

\def\Hommm#1#2#3{{\underline{\mathrm{Hom}}_{#1}(#2,#3)}}

\begin{document}
\selectlanguage{english} \catcode`\_=8\catcode`\^=7 \catcode`\_=8
\begin{abstract} We rewrite simplicially the standard definitions of
a complete first order theory, a model of it, and various characterisations of stability
of a complete first order theory. In our reformulations the simplicial language
replaces the standard definitions based on  syntax,
making them formally unnecessary. 

We view a complete first-order theory  as a symmetric simplicial object
in the category of profinite sets and open continuous maps,
defined by the functor sending a finite set of variables into the Stone space
of complete types in those variables. A model of a complete first-order theory is then
a morphism from a representable simplicial set satisfying certain lifting properties
reminiscent of, but weaker then, those in the definition of a fibration.
The class of simplicial profinite sets corresponding to complete first order theories
is characterised by the same lifting properties required of the map 
from the simplicial covering space (decalage) forgetting the extra degeneracy.
\end{abstract}
\maketitle

\def\todo{{\bf TODO:}}
\def\FinOrd{\text{FinOrders}}
\def\FiniteNonEmptySets{\text{FiniteNonEmptySets}^{\text{op}}} 
\def\op{\text{op}} 
\def\Aut{\text{Aut\,}}
\def\SS{\Bbb S} 
\def\SSS{\Bbb S}
\def\ACF{\text{ACF}} 
\def\ACFzero{\text{ACF}_0}
\def\ACFp{\text{ACF}_p}
\def\BbbC{\Bbb C}
\def\BbbR{\Bbb R}
\def\AutCQ{\Aut(\Bbb C/\Bbb Q)}
\def\AutRQ{\Aut(\Bbb R/\Bbb Q)}
\def\Vectk{{\text{Vect}_k}}
\def\GLVk{\text{GL}(\Vectk)}

\def\DopEmpty{\Delta^{\text{op}}_\emptyset}
\def\const{\text{const}} 

\def\ev{\text{ev}}
\def\Sym{\text{Sym}}
\def\AutL{\text{Aut}_{\mathcal L}}
\def\calL{\mathcal L}
\def\calM{\mathcal M}
\def\Bij{Bi\!j} 
\def\tp{\text{tp}}
\def\qftp{\text{qftp}}
\def\SSb{\Bbb S_\bullet} 

\def\St{{\ethi\ethmath{pa}}} 
\def\sSt{{s\!\St}}
\def\symSt{{sym\!\St}} 
\def\Top{\St}

\def\Stcl{{\ethi\ethmath{pO}}} 
\def\sStcl{{s\!\Stcl}}
\def\symStcl{{sym\!\Stcl}}

\def\Stclo{{\ethi\ethmath{pI}}} 
\def\sStclo{{s\!\Stclo}}
\def\symStclo{{sym\!\Stclo}}

\def\Stp{{\ethi\ethmath{pU}}} 
\def\sStp{{s\!\Stp}}
\def\symStp{{sym\!\Stp}}

\def\Stpo{{\ethi\ethmath{pE}}} 
\def\sStpo{{s\!\Stpo}}
\def\symStpo{{sym\!\Stpo}} 




%
%
\section{Introduction}
We transcribe the definitions of a complete first order theory, 
a model of a complete first order theory, 
and various equivalent characterisations of a stable first order theory,
in terms of the category of simplicial objects
of the category of profinite sets and open continuous maps, 
and related slightly larger categories. Simplicial language replaces
the syntax of the first order logic, making it formally unnecessary in this note. 
However, we do use the standard logic terminology {\em informally} to aid the intuition 
and clarify or state equivalence of our definitions and claims to 
the standard formulations in model theory. Finally, we need to say that the goal
of this short note is quite limited: to state the reformulations of a couple of 
basic definitions in model theory in most explicit simplicial language.
In particular, no attempt is made to relate these reformulations to 
the theory of simplicial sets or calculate cohomology.

By itself our reformulations lack any mathematical content 
and amount to little but a change of notation. However, they place these notions 
in a different context, which might eventually be useful to some. 
One potential application is that they may enable a more transparent, 
category-theoretic exposition of universal constructions such as used 
in R.Knight's  approach \cite{Knight} towards the Vaught conjecture. 
Another is that a more geometric understanding of the construction and uniqueness
of the pseudo-exponential function $\operatorname{ex}:\Bbb C^+\lra \Bbb C^*$ 
satisfying the Schanuel conjecture \cite[X\S6.15-18]{ManinZilber}, may eventually lead to 
a reformulation of the Grothendieck period conjecture, or 
something like pseudo-Betti cohomology satisfying some of the standard conjectures
on the motivic Galois group. 
One wonders how our simplicial reformulations of a first order theory
relate 
to those in terms of the cylindric algebras of Tarski, or of 
quasi-minimal structures and excellence \cite{Excellence}. 

Interpretations of complete first order theories are closely related to 
morphisms between them as simplicial profinite sets. A natural question
is then whether one can develop cohomological methods to prove 
that one first order theory does not interpret another, or, say, 
recover the groupoids or amalgamation properties 
defined in \cite[\S2,\S4]{HrushovskiGroupoids} or homology groups of types in \cite{GKK}.  

\subsection*{Acknowledgements} I thank Martin Bays, Assaf Hasson, Sharon Hollander, 
Will Johnson, Mark Kamsma, Robin Knight, Kobi Peterzil, Mikhail Rabinovich, 
Alex Usvyatsov, Boriz Zilber, 
for useful discussions. A discussion with Lukas Vokrinek lead 
to a reformulation of non-dividing in term of idempotents.  
This research was partially supported by the Israel   
Science Foundation grant No. 2196/20. I thank University of Haifa, and Kobi Peterzil, 
for hospitality. This research was partially supported by ISF grant 290/19.
\subsubsection*{A sketch of our reformulations} Our explanations remain self-contained 
if the reader ignores  the notions from logic and thus the next two sentences. 
With a first order theory $T$ we associate a functor $T_\bullet$ on non-empty finite sets 
which sends a finite set of variables into the Stone space of complete types in those variables.
With a  model $M\models T$ of $T$ we associate a map $|M|_\bullet\implies T_\bullet$ 
from the simplicial set 
$|M|_\bullet: n^\leqslant\longmapsto M^n$ represented by the of elements of $M$
which sends a tuple $(a_1,...,a_n)\in M^n$ into its type $\tp(a_1,...,a_n)$. 

To characterise the classes of functors arising in this way we first need 
to introduce a property of a map of simplicial sets 
equivalent to being fibrant for a map of fibrant simplicial sets, 
but in general weaker. Say that a map $X_\bullet \lra Y_\bullet$ is {\em vibrant} 
iff $X_\bullet\circ[-\infty\ll..] \lra Y_\bullet\circ[-\infty\ll..] \times_{Y_\bullet} X_\bullet$ 
is surjective where 
$[-\infty\ll..]:\Dop\lra\Dop, (1\!<..<\!n)\longmapsto (-\infty\!<\!1<..<\!n),n\geqslant 1$ denotes the 
decalage endomorphism adding the new least element to a finite linear order. 
A first order theory $T$ is then the same  as 
a symmetric simplicial profinite set $T_\bullet$ with open continuous simplicial maps 
such that the map $T_\bullet\circ[..\ll+\infty]\lra T_\bullet$ is vibrant.
An $<\!\aleph_0$-saturated model $M$ of $T$ is the same as a vibrant map 
$M_\bullet\lra T_\bullet$ from a representable simplicial set. 
Both these notions can also be defined by either left or right lifting properties 
with respect to explicitly given classes of morphisms. 
The various characterisations of stability of a first order theory are stated in 
terms of commutative diagrams in Theorem~\ref{thm:stable} and Theorems~\ref{thm:stableViaForking}.
One of these characterisation is a right lifting property.

\section{Simplicial definitions of f.o.theories}

In this section we shall transcribe the definitions of a first order theory, 
and of a model of a first order theory, 
in terms of the category of simplicial objects
of the category of profinite sets and open continuous maps, 
and related slightly larger categories. 

\subsection{Basic notation and conventions\label{subsec:basic_notation}}
We now introduce necessary simplicial notions and define 
the categories we use.  

\subsubsection{Simplicial notations} 
Let $n^\leqslant$ and $n-1$ denote a finite linear order with $n$ elements.
The categories $\Dop\subset \DopEmpty$ of non-empty (resp.~possibly empty) 
finite linear orders are often viewed as the full subcategories of the category of categories. 

We find useful the {\em decalage (shift)} endofunctors $[+1],[-1]:\Dop\lra \Dop $ 
where $[+1]=\const_0 \sqcup \id: (1<2<...<n) \longmapsto (-\infty<1<2<...<n)$ and
$[+1]=\id \sqcup \const_0 : (1<2<...<n) \longmapsto (1<2<...<n<\infty)$,
or rather the monads $[+1]\implies \id:\Dop\lra \Dop $ and $
[-1]\implies \id:\Dop\lra\Dop$. 
Sometimes we use the autofunctor $[-]:\Dop\lra\Dop, (1<2<...<n) \longmapsto (n<..<2<1)$ 
reversing the linear order. Evidently $[-]\circ[-]=\id$ and 
$[-1]=[-]\circ[+1]\circ[-]^{-1}$ and $[+1]$ are conjugated by $[-]$.
Sometimes we denote $[-\infty\ll..]:=[+1]$ and $[..\ll+\infty]:=[-1]$. 

\subsubsection{Notation for  categories of spaces}
Let $\Top\supset\Stp$ be the categories 
of topological, resp. profinite, spaces
with continuous maps.

Let $\Stpo$ be the category of profinite sets (=profinite spaces) with open continuous maps.
Let $\symStpo:=\Funct(\FiniteNonEmptySets , \Stpo)$ and
$\sStpo:=\Funct(\Dop , \Stpo)$ be its categories of symmetric simplicial and 
of simplicial objects. 

We also need a category to describe the structure inherited by a subset of a profinite space. 
A subset of a profinite space inherits the notion of a clopen subset 
which is closed under finite union and intersection, and generates
the induced topology on the subset. However, a clopen subset in the induced topology
does not necessarily come from a clopen subset of the profinite space:
for example, this often happens when the induced topology is discrete. 
Because of this, we define categories of {\em clopen pre-topological} spaces 
$\Stcl\supset\Stclo$ whose objects  are sets equipped with a collection of {\em clopen}
subsets closed under finite union and intersection. In $\Stcl$ morphisms are 
continuous maps (where the preimage of a clopen set is necessarily clopen), 
and in $\Stclo$ morphisms are continuous maps such that 
the closure of the image of a clopen subset is necessarily clopen.
We shall refer to such maps as {\em clopen continuous} maps. 

The above is summarised in the definition below. 
\begin{defi}[Categories of spaces]
Let $\Top\supset\Stp$ be the categories 
of topological, resp. profinite, spaces
with continuous maps.
Let $\Stcl\supset\Stclo$ 
denote the categories of {\em clopen pre-topological} spaces 
whose objects  are sets equipped with a collection of {\em clopen}
subsets closed under finite union and intersection. In $\Stcl$ morphisms are 
continuous maps, and in $\Stclo$ morphisms are continuous maps such that 
the closure of the image of a clopen subset is necessarily clopen.
We shall refer to such maps as {\em clopen continuous} maps. 
\end{defi}




\subsubsection{Vibrant simplicial maps} We shall define the notions of a first order theory and its $<\!\aleph_0$-saturated models in terms of  a property of maps of simplicial sets 
weaker than a fibration, and equivalent to being a fibration for a map of fibrant simplicial sets. 
In the category $\sFF$ and $\sFFae$ of generalised spaces 
a closely related property defines the notions 
of a map of topological spaces being closed or open. 

\begin{defi}[Virbrant simplicial maps] 
Call a map $M_\bullet\lra T_\bullet$ of simplicial sets {\em vibrant} iff 
it satisfies  either of the following equivalent conditions:
\bee
\item 
$ M_\bullet\circ[-\infty\ll ..\ll+\infty] \lra 
T_\bullet\circ[-\infty\ll ..] \times_{T_\bullet} M_\bullet\circ[..\ll+\infty]$
is surjective 
\item 
for any inclusion of finite sets $k\subset m,n$ 
the obvious map 
$$ M_\bullet((m+k+n)^\leqslant)   \lra 
T_\bullet({(m+k)^\leqslant})\times_{T_\bullet(k^\leqslant)} M_\bullet((k+n)^\leqslant)
$$ is surjective 
\item
$T_\bullet$ satisfies the following lifting property for each $0\leqslant k \leqslant m, n$ 
(below is the same equation  in two different notations)
\begin{equation}\label{TarskiVaughtRight}
\xymatrix@C=1.039in{ 
(k+n)^\leqslant_\bullet
 \ar[r] \ar[d]|\subseteq & M_\bullet \ar[d]|{(\therefore model)} \\
(m+k+n)^\leqslant_\bullet \ar@{-->}[ru] \ar[r] & T_\bullet}
\xymatrix@C=1.039in{ 
\Delta[k+n-1]_\bullet
 \ar[r] \ar[d]|\subseteq & M_\bullet \ar[d]|{(\therefore model)} \\
\Delta[m+k+n-1]_\bullet \ar@{-->}[ru] \ar[r] & T_\bullet }
\end{equation} 
\eee
\end{defi}
\begin{proof}
(1)$\iff$(2): surjectivity in dimension $k$ means exactly surjectivity of the map in (2) 
for $m=n=1$. Applying surjectivity of the latter map several times gives surjectivity for each 
finite $m,n\geqslant 0$.
(2)$\iff$(3): the diagonal arrows represents a simplex in the domain, 
and the commutative square represents a simplex in the product. 
\end{proof}

\begin{exem}[A fibration is vibrant] A fibration of simplicial sets is vibrant, and 
a map of fibrant simplicial sets is vibrant iff it is a fibration. 
This should be easy to see
by a diagram chasing calculation.
\end{exem}

\begin{exem}[Serre fibration as a vibrant map] A map $X\lra Y$ of topological spaces is a Serre fibration iff $\Singb X \lra \Singb Y$ is vibrant. 
More generally, a map $X_\bullet \lra Y_\bullet$ of fibrant simplicial sets is a fibration iff it is vibrant.    
\end{exem}
\begin{proof} This should be immediate to check. Surjectivity of 
$ \Singb X ((m+k+n)^\leqslant)   \lra 
 \Singb Y ({(m+k)^\leqslant})\times_{\Singb X(k^\leqslant)} \Singb X((k+n)^\leqslant)$
 for $m=1$ 
means  the lifting property of $X\lra Y$ wrt $\Delta[k+n] \lra \text{Cone}(\Delta[k+n])$
where $\Delta[k+n]\subset \Bbb R^{k+n+1}$ is the standard simplex as a topological space.
\end{proof}

\begin{exem}[Closed and open maps as being vibrant] 
The notions of a map of topological spaces being closed or open (=locally surjective) 
are captured by the notion of a vibrant map if we extend it to a category $\sFF$ of simplicial sets equipped 
with a notion of a ''simplex being small enough'', or, equivalently, a structure of topological nature giving a precise meaning to the phrase ``a property holds for all $n$-simplicies 
small enough''. 
We view $\sFF$ as a category of generalised spaces \cite{situs}, 
and discuss it in a forthcoming paper.

Formally, we consider simplicial objects in the category $\FF$ of sets equipped with 
a filter (= sets equipped with a finitely 
additive probabilistic measure taking only two values $0$ and $1$);
it is intuitive to refer to subsets of full measure as either {\em neighbourhoods} or {\em big subsets}. A morphism is a map such that the preimage of a neighbourhoods (=big subset) is 
a neighbourhood (=a big subset). In the definition of a vibrant map in $\sFF$ 
we require instead of surjectivity, so to say, local surjectivity: 
the image of a neighbourhood is a neighbourhood. The definition of a topological space
in terms of neighbourhoods leads to an embedding $\Topp\lra \sFF$ such that, in particular,
a map $X\lra Y$ is closed, resp. open, iff $X_\bullet\lra Y_\bullet$, resp. 
$X_\bullet\circ[-]\lra Y_\bullet\circ[-]$, is vibrant. 
\end{exem}

\subsection{A first order theory and its models, simplicially\label{subsec:basic_defs}}
We shall now transcribe simplicially basic definitions in first order logic
as objects in simplicial categories $\Stpo\subset \Stclo$ satisfying certain lifting 
properties. Recall that examples of these definitions are given in the 
next subsection \S\ref{subsec:basic_examples}. 

\subsubsection{Introduction} 
We view a first order theory $T$
as the simplicial quotient 
$$T_\bullet:= \faktour{|M|_\bullet}{\Aut(M)}$$ $$
T_\bullet(n^\leqslant):= \faktour{M^n}{\Aut(M)} \text{(by the diagonal action)}$$ 
of a sufficiently large and homogeneous ``monster'' model $M$ of $T$.
This data packs together  the  spaces of {\em complete $T$-types} with Stone topology
where maps in $\Dop$ correspond to manipulations of the variables.\footnote{%
Let $|M|_\bullet$, 
$|M|_\bullet(n^\leqslant):=\Homm {Sets} {n^\leqslant} $ denote 
the simplicial set 
represented by the set of points of $M$. For a group $G$ acting on a set $M$, 
let $|M|_\bullet/G$ denote {\em the simplicial quotient}
$\faktour{|M|_\bullet}G (n^\leqslant):=\faktour{|M|_{n-1}}G=
\faktour{M^{n}}G,\,\, n\geqslant 1$.
For example, when $M:=G$ is $G$ acting on itself by multiplication, 
this quotient is the simplicial classifying space of a group $G$  
$\BBBb G:= \faktour{|G|_\bullet}G$.} 
The first order formulas of $T$, or rather their equivalence classes defining the same subset, 
correspond to the clopen subsets of $T_\bullet(n^\leqslant)$, $n\geqslant 0$. 

This data remembers everything about the first order theory but which formulas correspond to the basic predicates in its language. 
We give several examples in \S\ref{subsec:basic_examples}.

This is very similar to the notions of 
type structure \cite{Morley71}, 
 type category
 \cite{Knight}[Definition 2.6 and Theorem 2.8] 
 and type space functor \cite{Kamsma}[Def.4.13].
%

\subsubsection{A first order theory in the form of Knight} A first order theory is 
a symmetric simplicial object $T_\bullet:\Dop\lra\Stpo$ 
in the category of clopen maps of the profinite sets
such that its decalage $T_\bullet[..\ll+\infty]\lra T_\bullet$ is vibrant. 
As we shall see later, a $<\!\aleph_0$-saturated model of a f.o.theory 
can be defined as 
a vibrant map $M_\bullet \lra T_\bullet$ from a representable symmetric simplicial set.
Thus, so to say, a first order theory is
a symmetric simplicial object in $\symStclo$ such that its shift 
is its $<\!\aleph_0$-saturated model.

\begin{defi}[f.o.theory]\label{defi:fotheory} A {\em first order theory} $T$ 
is the same as a symmetric simplicial 
object of the category $\Stpo$ of profinite sets and open(=clopen) continuous  maps, 
i.e. a functor 
$$T_\bullet:\FiniteNonEmptySets \lra \Stpo, $$
such that the map $T_\bullet\circ[..\ll+\infty]\lra T_\bullet$ is vibrant. 

That is, if it satisfies
 either of the following equivalent conditions:
\bee
\item 
$ T_\bullet\circ[-\infty\ll ..\ll+\infty] \lra 
T_\bullet\circ[-\infty\ll ..] \times_{T_\bullet} T_\bullet\circ[..\ll+\infty]$
is surjective 
\item 
for any inclusion of finite sets $k\subset m,n$ 
the map 
$$ T_\bullet((m+k+n)^\leqslant) \lra 
T_\bullet({(m+k)^\leqslant})\times_{T_\bullet(k^\leqslant)} T_\bullet((k+n)^\leqslant)
$$ is surjective 

\item the map $T_\bullet\circ[..\ll+\infty]\lra T_\bullet$ satisfies the following lifting property for each $0\leqslant k \leqslant m, n$ 
(below is the same equation  in two different notations)
\begin{equation}\label{TarskiVaughtRight}
\xymatrix@C=1.039in{ 
(k+n)^\leqslant_\bullet
 \ar[r] \ar[d]|\subseteq & T_\bullet\circ[+m] \ar[d] \\
(m+k+n)^\leqslant_\bullet \ar@{-->}[ru] \ar[r] & T_\bullet}
\xymatrix@C=1.039in{ 
\Delta[k+n-1]_\bullet
 \ar[r] \ar[d]|\subseteq & T_\bullet\circ[+m] \ar[d] \\
\Delta[m+k+n-1]_\bullet \ar@{-->}[ru] \ar[r] & T_\bullet }
\end{equation} 

\item
$T_\bullet$ satisfies the following lifting property for each $0\leqslant k \leqslant m, n$ 
(below is the same equation  in two different notations)
\begin{equation}\label{TarskiVaughtRight}
\xymatrix@C=1.039in{ 
(m+k)^\leqslant_\bullet\vee_{k^\leqslant_\bullet} (k+n)^\leqslant_\bullet
 \ar[r] \ar[d]|\subseteq & T_\bullet \\
(m+k+n)^\leqslant_\bullet \ar@{-->}[ru] & }
\xymatrix@C=1.039in{ 
\Delta[m+k-1]_\bullet\vee_{\Delta[k-1]_\bullet} \Delta[k+n-1]_\bullet
 \ar[r] \ar[d]|\subseteq & T_\bullet \\
\Delta[m+k+n-1]_\bullet \ar@{-->}[ru] & }
\end{equation} 
\eee

A {\em first order theory in a countable language}  
is a symmetric simplicial 
object in the category of light (=metrisable) profinite sets and 
open(=clopen) continuous maps.

A {\em formula with $n$ variables} $v_1,...,v_n$ is the same as 
an clopen subset of $T_\bullet(\{v_1,...,v_n\})$.
Logic connectives {\em conjunction $\varphi(v_1,...,v_n) \aand \psi(v_1,...,v_n)$, 
disjunction $\varphi(v_1,...,v_n) \vee \psi(v_1,...,v_n)$}, and 
{\em negation $\neg\varphi(v_1,...,v_n)$},
correspond to {\em intersection, union,} and {\em complement} of the clopen subsets
of $T_\bullet(\{v_1,...,v_n\})$. 
The {\em existential quantifier $\exists v_i \varphi(v_1,...,v_n)$} 
correspond to taking the image under the inclusion map 
$\{v_1,..,v_{i-1},v_{i+1},...,v_n\}\subseteq 
\{v_1,..,v_{i-1},v_i,v_{i+1},...,v_n\}$. That is, 
\begin{equation}\begin{split}
 \varphi(v_1,...,v_n) \aand \psi(v_1,...,v_n)  &:=
 \varphi(v_1,...,v_n) \cap \psi(v_1,...,v_n) \\ 
 \varphi(v_1,...,v_n) \vee \psi(v_1,...,v_n)  &:=
\varphi(v_1,...,v_n) \cup \psi(v_1,...,v_n) \\
 \neg\varphi(v_1,...,v_n)  &:=  T_\bullet(\{v_1,...,v_n\})\setminus \varphi(v_1,...,v_n) \\
 \exists v_i \varphi(v_1,...,v_n)  &:=  \Imm \left(
T_\bullet( \{v_1,..,v_{i-1},v_i,v_{i+1},...,v_n\})
\xra {} 
T_\bullet(\{v_1,..,v_{i-1},v_{i+1},...,v_n\})
\right)\end{split}\end{equation}
\end{defi}
\begin{proof} In set theoretic language, this is  
\cite{Knight}[Definition 2.6 and Theorem 2.8]
The lifting property \eqref{TarskiVaughtRight} appears under the name {\em amalgamative} 
in Definition 2.6, p.58. In the standard logic terminology  \eqref{TarskiVaughtRight} says 
that for any two $T$-types
$p(x_1,...,x_m,y_1,...,y_k)$ and $q(y_1,...,y_k,z_1,...,z_n)$ there is a type 
$$r(x_1,...,x_m,y_1,...,y_k,z_1,....,z_n)\supset 
p(x_1,...,x_m,y_1,...,y_k), q(y_1,...,y_k,z_1,...,z_n)$$ 
extending them both. In other words, if you have tuples $a,b$ and $b,c$ and 
thus their types $\tp(a,b)$ and $\tp(b,c)$, then you may form a tuple $a,b,c$
and thereby its type $\tp(a,b,c)$. \end{proof}

%
\begin{exem}[A nerve and a classifying space of a profinite group] 
A nerve $N_\bullet$ of a category is 
an example of a non-symmetric simplicial set 
which satisfies the defining property of a f.o.theory 
that $N_\bullet\circ[..\ll+\infty]\lra N_\bullet$ is surjective.
So is the singular simplicial set $\Singb X$ of a topological space, 
and, more generally, any fibrant simplicial set. 

The simplicial model $\BBBb G:= \faktour{|G|_\bullet}G$ of the classifying space
of a profinite group $G$ is an example of a profinite simplicial set
satisfying the defining property of a f.o.theory (being a nerve) but 
the diagonal maps are usually not open.
%
\end{exem}

\begin{que} Find a natural examples of a non-symmetric ``f.o.theory'', i.e. 
a simplicial profinite set with clopen simplicial maps such that its shift is vibrant.
\end{que}

\subsubsection{A $<\!\aleph_0$-saturated model of a f.o.theory in the form of Knight} 
Recall that a model $M$ of a first order theory $T$ is {\em $<\!\aleph$-saturated}
iff for each tuple $b\in M^n$ each type $p(x,b)$ is realised in $M$, i.e. there is 
a tuple $a\in M^m$ such that $p(x,b)$ is the type of the tuple $a,b$. 
This notion is captured by the definition below. 
\begin{defi}[a $<\!\aleph_0$-saturated model of f.o.theory]\label{defi:fotheory} 
A {\em $<\!\aleph_0$-saturated model $M$} of a first order theory $T$ 
is the same as a vibrant map to $T_\bullet$ 
from a 
representable symmetric simplicial object of the category $\Stclo$ 
of subsets of profinite sets and clopen continuous  maps.

That is, it is a natural transformation 
$$M_\bullet\implies T_\bullet:\FiniteNonEmptySets \lra \Stpo, $$
satisfying either of the equivalent conditions:
\bee
\item 
$ M_\bullet\circ[-\infty\ll ..\ll+\infty] \lra 
T_\bullet\circ[-\infty\ll ..] \times_{T_\bullet} M_\bullet\circ[..\ll+\infty]$
is surjective 
\item 
for any inclusion of finite sets $k\subset m,n$ 
the map 
$$ M_\bullet((m+k+n)^\leqslant) \lra 
T_\bullet({(m+k)^\leqslant})\times_{T_\bullet(k^\leqslant)} M_\bullet((k+n)^\leqslant)
$$ is surjective 
\item
$T_\bullet$ satisfies the following lifting property for each $0\leqslant k \leqslant m, n$ 
(below is the same equation  in two different notations)
\begin{equation}\label{TarskiVaughtRight}
\xymatrix@C=1.039in{ 
(k+n)^\leqslant_\bullet
 \ar[r] \ar[d]|\subseteq & M_\bullet \ar[d]|{(\therefore model)} \\
(m+k+n)^\leqslant_\bullet \ar@{-->}[ru] \ar[r] & T_\bullet}
\xymatrix@C=1.039in{ 
\Delta[k+n-1]_\bullet
 \ar[r] \ar[d]|\subseteq & M_\bullet \ar[d]|{(\therefore model)} \\
\Delta[m+k+n-1]_\bullet \ar@{-->}[ru] \ar[r] & T_\bullet }
\end{equation} 
\eee
\end{defi}

\subsubsection{A f.o.theory in the form of Kamsma} 
The definition below describes the class of f.o.theories by a lifting property
on the other side. It aims to describe the configurations of open subsets 
forbidden by the Beck-Chevalley condition. In terms of the previous definition,
these open subsets are clopen subsets containing the simplicies specified 
by the maps in \eqref{TarskiVaughtRight}, namely $y_{m+k+n}\in Y_{(m+k+n)^\leqslant}$
in \eqref{TarskiVaughtLeftBC} below would be the simplex
corresponding to the lifting arrow in \eqref{TarskiVaughtRight}. 
The notation $(Y_\bullet: Y_{k^\leqslant}:=Y_{k^\leqslant}\setminus \{y_{k}\})$
denotes the simplicial set $Y_\bullet$ with the simplex $y_{k}\in Y_{k^\leqslant}$
punctured out, along with all its preimages. 
\begin{defi}[f.o.theory]\label{defi:fotheory} A {\em first order theory} $T$ 
is the same as a symmetric simplicial 
object of the category $\Stpo$ of profinite sets and open(=clopen) continuous  maps, 
i.e. a functor 
$$T_\bullet:\FiniteNonEmptySets \lra \Stpo$$ 
whose underlying simplicial set is connected, and 
such that
for each connected finite discrete symmetric simplicial set $Y_\bullet$, 
each (not necessarily commutative) square in $\FiniteNonEmptySets$ 
as shown on the left, 
and for each pair of points $y_{m+k}\in Y_\bullet((m+k)^\leqslant)$ and $
y_{k+n}\in Y_\bullet((k+n)^\leqslant)$ such that
$T(\rho_{m+k\rightarrow k})(y_{m+k})=T(\rho_{k+n\rightarrow k})(y_{k+n})$
the following lifting property  holds
\begin{equation}\label{TarskiVaughtLeftBC}
\xymatrix@C=1.8239in{
(m+k+n)^\leqslant \ar@{<-}[r]|>>>>>>>{\rho_{k+n\rightarrow  m+k+n}} \ar@{<-}[d]_{\rho_{m+k\rightarrow m+k+n}}
^{\,\,\,\,\rho_{k\rightarrow m+k}(y_{m+k})=\rho_{k\rightarrow k+n}(y_{k+n})}
& 
(k+n)^\leqslant \ar@{<-}[d]^{\rho_{k\rightarrow  k+n}}
\\ (m+k)^\leqslant \ar@{<-}[r]|{\rho_{k\rightarrow  m+k}} & k^\leqslant } 
\end{equation}
implies

			  \begin{equation}
			  \begin{tikzcd}[row sep=large, column sep=tiny]
			    & (Y_\bullet: Y_{(m+k)^\leqslant}:=Y_{(m+k)^\leqslant}\setminus \{y_{m+k}\}) \sqcup (Y_\bullet: Y_{(k+n)^\leqslant}:=Y_{(k+n)^\leqslant}\setminus \{y_{k+n}\}) \arrow[d] \\
			      T_\bullet \arrow[r] \arrow[ru, dashed] & 
			        (Y_\bullet: Y_{(m+k+n)^\leqslant}:=Y_{(m+k+n)^\leqslant}\setminus (\rho_{m+k\rightarrow m+k+n}^{-1}(y_{m+k})\cap \rho_{k+n\rightarrow m+k+n}^{-1}(y_{k+n})))
				\end{tikzcd}
				\end{equation}
%
%
\end{defi}
\begin{proof} 
For simplicial profinite sets 
the lifting property \eqref{TarskiVaughtLeftBC} should be equivalent to 
 \eqref{TarskiVaughtRight} and thereby to
the Beck-Chevalley condition as stated in \cite{Kamsma}.
This should be equivalent \cite{Kamsma}[Def. 4.7, Fact 4.8].
\end{proof}


\subsubsection{A model of a f.o.theory} A model of a f.o.theory 
is a map from a representable set to the f.o.theory satisfying 
a certain lifting property expressing something like surjectivity. 
However, that representable set no longer 
lives in the category $\symStpo$ as its objects are no longer profinite.  
\begin{defi}[model]\label{defi:model} 
A {\em model} of a f.o.theory $T_\bullet$ is a morphism $M_\bullet\lra T_\bullet$
in the category of clopen pre-topological  spaces $\symStclo$ with clopen continuous maps,
from a symmetric simplicial space $M_\bullet$ 
whose underlying simplicial set $|M|_\bullet$ is representable, and satisfying the 
lifting property $\eqref{TarskiVaughtLeft}$ 
\begin{equation}\label{TarskiVaughtLeft}\begin{split} 
\xymatrix@C=2.39in{
|M|_\bullet\circ[+n]\ar[r]\ar[d] & 
|\{\bullet\!\!\!\rlap{\raisebox{2.78ex}{$|\{\star_2\}|_\bullet$}}\,\,\,,\star_1\}|_\bullet 
\ar[d] \\
|T|_\bullet\circ[+n]\ar@{-->}[ru]\ar[r] & |\{\bullet,\star\}|_\bullet} 
\\ \text{\tiny the topology is such that simplicies }(\star_1,...,\star_1)
\text{ \tiny and }(\star_2,...,\star_2)\text{ \tiny are indistinguishable}
\end{split}\end{equation}
\end{defi}

Recall that an {\em $<\!\!\aleph_0$-saturated model} is a morphism $|M|_\bullet\lra T_\bullet$ 
in $\symStclo$ from a representable symmetric simplicial space 
$|M|_\bullet$ satisfying 
the same lifting property $\eqref{TarskiVaughtRight}$ appearing in the definition 
of a f.o.theory.  
\begin{equation}\label{model_fibration}\xymatrix{
M_\bullet\circ[+1] \ar[ddr] \ar@{->>}[dr]|{(surjective)} \ar[drr] & \\
&T_\bullet\circ[+1] \times_{T_\bullet} M_\bullet \ar[d] \ar[r] & M_\bullet \ar[d] & \\
&T_\bullet[+1] \ar[r] & T_\bullet &
}\end{equation} 

\begin{proof} 
Consider \eqref{TarskiVaughtLeft}.
As a simplicial set 
$|M|_\bullet\circ[+n]=\bigsqcup\limits_{(a_1,...,a_n)\in M^n} 
\{(a_1,...,a_n)\}\times |M|_\bullet$ is the disjoint union of 
simplicial connected components,
and so is 
$T_\bullet\circ[+n]=\bigsqcup\limits_{p_n\in T_\bullet(n^\leqslant)} 
\{q_\bullet\in T_\bullet:q_\bullet[1,...,n]=p_n\}$.
On a connected component of $T_\bullet\circ[+n]$ with empty preimage,
define the lifting trivially by sending everything to the connected component of 
$|\{\bullet,\star_1\}|_\bullet $
of $|\{\bullet,\star_1\}|_\bullet \bigsqcup |\{\star_2\}|_\bullet$.
Therefore we only need to consider the connected components of $T_\bullet(n^\leqslant)$
with non-empty preimage in $|M|_\bullet\circ[+n]$, i.e. 
where the type $p_n(v_1,...,v_n)$ is realised in $M$.
The preimage of $T_\bullet((n+1)^\leqslant) \lra 
|\{\bullet,\star\}|_\bullet $
of $ |\{\star\}|_\bullet$ is an arbitrary clopen subset, i.e. 
a formula $\varphi(x_1,...,x_n,y)$ of $T$. 
There is a commutative square where 
the whole fibre $\{(a_1,...,a_n)\}\times |M|_\bullet$ is sent to $\{\star_2\}$ 
iff empty is the preimage 
of $ |\{\bullet\}|_\bullet$
under the map $\{(a_1,...,a_n)\}\times |M|_\bullet\lra 
|\{\bullet,\star\}|_\bullet $.
That is, iff the formula $\varphi(a_1,...,a_n,y)$ is not realised in $M$. 
There is a lifting 
iff $\tp(a_1,...,a_n)\aand \varphi(a_1,...,a_n,y)$ is inconsistent 
in $T$, i.e. $\neg \exists y \varphi(a_1,...,a_n,y)\in \tp(a_1,...,a_n)$. 
In other words, we have the implication 
$ \exists y \varphi(a_1,...,a_n,y)\in \tp(a_1,...,a_n)$
implies $\varphi(a_1,...,a_n,y)$ is realised in $M$.
This should be equivalent to the Tarski-Vaught condition for $M$. 
Recall a structure $M$ in the language $L$ of  a f.o.theory $T$
is a model of $T$ iff $M$ satisfies the {\em Tarski-Vaught condition} 
saying that 
for each formula $\varphi(\bar x,\bar y)$ of theory $T$ and 
for each tuple $\bar a \in M^n$  
\bi\item[\eqref{TarskiVaughtLeft}$_\text{words}$] 
if $M\models\exists \bar y \varphi(\bar a,\bar y)$ 
then for some tuple $\bar b\in M$ $M\models \varphi(\bar a,\bar b)$ 
\item[] or equivalently
if $\exists y \varphi(x,y)\in \tp(a)$ 
then for some tuple $\bar b\in M$ $M\models \varphi(\bar a,\bar b)$ 
\ei

Hence, the lifting property \eqref{TarskiVaughtLeft} should be equivalent 
to the Tarski-Vaught condition \eqref{TarskiVaughtLeft}$_\text{words}$. 

The lifting property \eqref{TarskiVaughtRight} should be equivalent 
to the Tarski-Vaught condition for types rather than formulas. In fact, 
it should be equivalent if we assume $m=0$. 
Recall that the standard Tarski-Vaught criterion for types says:
if $p(x,y)$ is a type, and $a$ is a tuple in $N$, then there is a tuple $b$ in $N$
realising type $p(a,y)$. As a lifting property this is 
for any $0\leqslant k \leqslant n$ it holds 
\begin{equation}\label{TarskiVaughtForTypesAndModels} 
\xymatrix@C=2.39in{  k^\leqslant_\bullet \ar[r] \ar[d]|\subseteq & |N|_\bullet \ar[d] \\
 (k+n)^\leqslant_\bullet \ar@{-->}[ru] \ar[r] & T_\bullet }\end{equation} 
\end{proof}

\subsection{A f.o.theory as a profinite completion of its model}
We sketch a simplicial definition of a structure in a first order language,
and then ask whether there is a precise sense in which 
the f.o.theory of a structure is its profinite completion/compactification. 

\subsubsection{A structure in language $\mathcal L$, simplicially} 
Recall that {\em a first-order structure $M$ in a language $\mathcal L$}
is a set $M$ equipped with several distinguished subsets of 
finite Cartesian powers of $M$. These distinguished subsets are 
called {\em predicates} of the language $\mathcal L$. 
Take the simplicial set represented by the set $|M|$ of elements of $M$, 
and equip each $|M|_\bullet(n^\leqslant):=|M|^n$ with the 
coarsest ``clopen pre-topology", i.e. 
a collection of ``clopen'' subsets closed under finite union and intersection, 
such that these distinguished subsets are clopen, and so are all the symmetric simplicial maps
$|M|^n\lra |M|^m$. 
The clopen subsets correspond precisely to the first order formulas in $\mathcal L$
without parameters, or, equivalently, to subsets first order definable in $\mathcal L$ 
without parameters.
We shall denote
by $M^{\mathcal L}_\bullet$ or simply $M_\bullet$ the resulting object of 
the category $\symStclo$ 
of symmetric simplicial clopen pre-topological spaces with continuous maps
such that the closure of the image of a clopen subset is necessarily clopen.

If $M$ is a saturated model of its first order theory in the language $\mathcal L$, 
this is the same as equipping each $|M|_\bullet(n^\leqslant):=|M|^n$ with the 
coarsest topology with the same properties. 

However, if $M$ is not saturated, this may fail. It still holds that 
any first order $\emptyset$-definable subset is necessarily clopen, 
but the converse is no longer true. For example, for
$M:=(\Bbb N,\,+\subset \Bbb N\times \Bbb N\times \Bbb N)$ the additive group of natural numbers,
each natural number is definable without parameters, and therefore each point is clopen.
Hence, the resulting topology is discrete and thus remembers no information about the language.



\subsubsection{The f.o.theory of a $<\!\aleph_0$-saturated model is 
its profinite compactification ?} 
There is a coaugmented compactification functor $\St\lra \Stp$
from the category of topological spaces into the category of profinite spaces
$$ X\longmapsto \cl \Imm \left( X \lra \prod_{X\lra\{0,1\}} \{0,1\} \right) $$
fitting into a weak factorisation system \cite{QuillenNegation}[Fig.5,p.18]
\begin{equation} X \xra {\left\{\mathcalAnobrackets\,\,\,\nondenseimage\right\}^l} 
X^{\Stp}
\xra{\left\{\mathcalAnobrackets\,\,\,\nondenseimage\right\}^{lr}} 
\bot \end{equation}
This functor should induce  a compactification functor 
$\symStclo\lra\symStpo$ from the category of clopen pre-topological spaces 
to that of profinite spaces, both with clopen continuous maps.
A verification should show  that this profinitisation functor takes 
a structure $M^{\mathcal L}$ into its f.o.theory $T^{M^{\mathcal L}}_\bullet$,
and gives the coaugmentation map 
\begin{equation} 
M^{\mathcal L}_\bullet \lra T^{M^{\mathcal L}}_\bullet\end{equation}
which should fit into an weak factorisation system, say generated by 
the morphisms used in the lifting property in the definition of a model
or a f.o.theory.

\subsection{Interpretations of f.o.theories as their morphisms}

Our simplicial definition of a f.o.theory gives rise to a notion of 
{\em morphisms of f.o.theories} which roughly corresponds
to an interpretation of one f.o.theory in terms of another f.o.theory. 
This relationship is captured by a notion of 2-category of categories 
\cite[Def.3.6(2-category of theories),p.6]{Kamsma} in a more general context.

Below $M_\bullet\lra T_\bullet$ denotes a $<\!\aleph_0$-saturated model of a f.o.theory $T$.

\subsubsection{A definable subset with induced structure as a morphism} 
Let $D:=\{a \in M \,:\, M\models\phi(a)\}=\phi(M)$ be a definable subset of $M$.
Equip $D$ with the structure inherited from $M$, i.e. 
equip $D_\bullet\subset M_\bullet$ with the notion of clopen subset (=f.o.formula)
inherited from $M$. By definition this gives a morphism $D_\bullet\lra M_\bullet$ 
in $\sStclo$, and thereby a $\sStclo$-morphism $T^D_\bullet\lra T_\bullet$
of the f.o.theories where  $p(x) \longmapsto p(x) \& \phi(x)$.

Note that there is a partial map 
$T_\bullet  \dashrightarrow T^D_\bullet$ undefined on types $p(x)\ni\neg\phi(x)$.

Let $\phi(x,y)$ be a formula defining an equivalence relation on $M$. 
It induces a map $M \lra \faktour{M}{\approx_\phi}$. Equip $\faktour{M}{\approx_\phi}$
with the structure induced from $M$, i.e. consider the quotient map 
$M_\bullet \lra (\faktour{M}{\approx_\phi})_\bullet$ in $\sStclo$. It induces
a map of f.o.theories $T_\bullet \lra T^{\phi(-,-)}_\bullet$ in $\sStclo$.

\subsubsection{A reduct $\mathcal L'\subset \mathcal L$  as a morphism}
Let $T^{\mathcal L'}$ be the f.o.theory consisting of formulae in a sublanguage $\mathcal L' \subset \mathcal L$ of the language $\mathcal L$ of f.o.theory $T$. 
Evidently the inclusion ${\mathcal L'}\subset \mathcal L$ induces  a map 
$T_\bullet \lra T_\bullet^{\mathcal L'}$.

\subsubsection{A morphism of f.o.theories is a conservative extension of a 
definable subset ?} 
It seems that a morphism is something like a conservative extension of a subset, so 
in particular a pure field in a theory $Q$ is something like 
a map $(ACF_q)_\bullet\lra Q_\bullet$. Let me try to explain this. 

Consider a morphism $\xi_\bullet:T_\bullet\lra Q_\bullet$ in $\symStpo$. 
(i) By continuity each formula (=clopen subset) 
$\varphi(v_1,...,v_n)\subseteq Q_\bullet(\{v_1,...,v_n\})$ of f.o.theory $Q$ 
defines a formula (=the preimage of the clopen subset) 
$\xi\inv\varphi(v_1,...,v_n)\subseteq T_\bullet(\{v_1,...,v_n\})$
of f.o.theory $T$. 

(ii) In $\symStpo$ we also require the maps to be open, and that means that
each formula 
$\psi(v_1,...,v_n)\subseteq T_\bullet(\{v_1,...,v_n\})$
of f.o.theory $T$
defines a formula (=the image of the clopen subset) 
$\xi\psi(v_1,...,v_n)\subseteq Q_\bullet(\{v_1,...,v_n\})$ of f.o.theory $Q$.
In other words, a $T$-definable subset is also necessarily $Q$-definable. 

In the standard model theoretic terminology   
(i) should  means $\xi\inv$ defines {\em an interpretation of a definable subset
$\xi(x=_T x)$ of f.o.theory $Q$}.  
And (ii) should mean that each $T$-definable subset is also $Q$-definable, and 
thereby this interpretation should 
define a {\em conservative extension of a definable subset  of $Q$} in 1 variable.  

\section{Basic examples of f.o.theories\label{subsec:basic_examples}} 
We now give a couple of standard examples illustrating the definitions above. 
Namely, we define the f.o.theory $\ACFzero$ of an algebraically closed field,
of the real closed field $\BbbR$, and of the real field extended by the exponential function and all the real analytic functions restricted to an interval. 
\subsection{Simplicial quotient by a group action} 
Let $|M|_\bullet$, 
$|M|_\bullet(n^\leqslant):=\Homm {Sets} {n^\leqslant} $ denote 
the simplicial set 
represented by a set $|M|$ of points of some space or structure $M$. 
For a group $G$ acting on a set $|M|$, 
let $|M|_\bullet/G$ denote {\em the simplicial quotient}
$\faktour{|M|_\bullet}G (n^\leqslant):=\faktour{|M|_{n-1}}G=
\faktour{M^{n}}G,\,\, n\geqslant 1$.
For example, when $M:=G$ is a group $G$ acting on itself by multiplication, 
this quotient is the simplicial classifying space of a group $G$  
$\BBBb G:= \faktour{|G|_\bullet}G$.

If a group $G$ acts on $M$ by automorphisms, then 
sometimes the simplicial quotient inherits some of the structure of $M$. 
In that case we shall often write $\faktour{M_\bullet}{G}$ after making 
an appropriate definition. 

\subsection{The f.o.theory $\ACFzero$ of an algebraically closed field}
Let $\BbbC$ denote the field of complex numbers, and let $\AutCQ$ be 
its field automorphism group. 
As a first order structure, we consider the complex field $\BbbC$ in the language of fields
$$\mathcal L:= \big\{ \{(x,y,z):x+y=z\},\, \{(x,y,z):xy=z\} \big\}$$
consisting of two distinguished subsets (=predicates). 
Recall that the $\mathcal L$-structure is the 
object 
$\BbbC_\bullet^{\mathcal L} $ of $\symStclo$ whose underlying simplicial set is $|\BbbC|_\bullet$ is represented by the set of points of the field $\BbbC$, and 
the topology on each $\BbbC^n$ is the coarsest topology such that
these two subsets are both closed and open. It is known 
that in this case the clopen subsets of $\BbbC^n$ are 
precisely the constructible subsets defined over $\BbbQ$.

The {\em first order theory of an algebraically closed field of characteristic 0} 
is then the simplicial quotient 
\begin{equation}\faktour{\BbbC_\bullet}{\AutCQ}\end{equation}
where the clopen subsets of $\faktour{\BbbC_\bullet}{\AutCQ}(n^\leqslant)=
\faktour{\BbbC^n}{\AutCQ}$ are precisely the quotients of 
the constructible subsets of $\BbbC^n$ defined over $\BbbQ$. 
An orbit of the diagonal action of $\AutCQ$ on $\BbbC^n$ 
consists of the generic points of a Zariski-closed subvariety of $\BbbC^n$ 
irreducible over $\BbbQ$. Thus, one may think that 
{\em an $n-1$-simplex of the f.o.theory $\BbbC_\bullet$  
is an irreducible Zariski closed closed subvariety of $\BbbC^n$}.
In particular, the number of simplicies is countable, and thus
the f.o.theory $\faktour{\BbbC_\bullet}{\AutCQ}$ is a light profinite set.

\subsection{The f.o.theory $\BbbR$ of the real closed field}
Let $\BbbR$ denote the field of real numbers, and let $\AutRQ$ be 
its field automorphism group. 
As a first order structure, we consider the real field $\BbbR$ in the language 
$$\mathcal L:= 
\big\{ \{(x,y,z):x+y=z\},\, \{(x,y,z):xy=z\},\, \{(x,y):x\leqslant y\}\, \big\}$$
consisting of $3$ distinguished subsets (=predicates). 
Recall that the $\mathcal L$-structure is the 
object 
$\BbbR_\bullet^{\mathcal L} $ of $\symStclo$ whose underlying simplicial set 
is $|\BbbR|_\bullet$ is represented by the set of points of the field $\BbbR$, 
and the topology on each $\BbbR^n$ is the coarsest topology such that
these two subsets are both closed and open. It is known 
that in this case the clopen subsets of $\BbbR^n$ are 
precisely the Boolean combinations of real algebraic subsets of $\BbbR^n$
defined over $\Bbb Q$.

The {\em first order theory of a real closed field of characteristic $0$} 
is then the profinite completion of 
$\BbbR^{\mathcal L}_\bullet$ where clopen subsets are precisely
the Boolean combinations of real algebraic subsets of $\BbbR^n$ 
defined over $\BbbQ$. 
In particular, this profinite compactification adds non-standard, 
infinitely small 
or large elements, or rather their types, as the intersection
of the infinite families 
$\{x:0<x<\varepsilon\}_{\varepsilon>0}$
and $\{x:x>r\}_{r\in \BbbR}$ of clopen subsets.

\subsection{Stable theories: real vs complex numbers} 
Real and complex numbers provide two important examples of a unstable and of stable theory, i.e. 
the profinite completion $\BbbR^{\text{field}}_\bullet$ 
fails the equivalent conditions of 
Theorem~\ref{thm:stable} whereas $\faktour{\BbbC_\bullet}{\AutCQ}$ satisfies them.

\subsection{The f.o.theory $\BbbR^{\text{exp,an}}$ of the real field with analytic functions}
We still get a well-studied structure \cite{vdDries1998}
if we significantly extend the language 
and add (the graphs of) the real exponential function $\exp:\BbbR\lra\BbbR$
and all the real analytic functions $[a,b]\lra \BbbR$, $a,b\in \BbbR$ 
restricted to a finite closed interval. A key feature is that
a clopen subset (= definable with parameters) of each connected component of 
$\BbbR^{\text{exp,an}}_\bullet\circ[+n](1^\leqslant)$ is a finite union of intervals (= definable 
with parameters in $(\BbbR,<)$). The f.o.theories with this property are called
{\em o-minimal}.

\subsubsection{The f.o.theory $\Vectk$ of a vector space over a field $k$}
Let $\Vectk$ denote a vector space over the field $k$ of infinite dimension, 
and let $\GLVk$ be 
its automorphism group. 
The {\em first order theory of the vector space over a field $k$} 
is then the simplicial quotient
\begin{equation}\faktour{\Vectk_\bullet}{\GLVk}\end{equation}
where the clopen subsets of $\faktour{\Vectk_\bullet}{\GLVk}(n^\leqslant)=
\faktour{\Vectk^n}{\GLVk}$ are precisely  the quotients of 
Boolean combinations of $k$-affine subspaces  of $\Vectk^n$.

As a first order structure, we consider the field $\Vectk$ in the language 
$$\mathcal L:= \big\{ \{(x,y,z):x+y=z\} \big\}\cup  \{(x,y):y=ax\}_{a\in k}$$
Recall that the $\mathcal L$-structure is the 
object 
$\Vectk_\bullet^{\mathcal L} $ of $\symStclo$ whose underlying simplicial set is 
$|\Vectk|_\bullet$ is represented by the set of points of the vector space $\Vectk$, 
and the topology on each $\Vectk^n$ is the coarsest topology such that
these two subsets are both closed and open. It is known 
that in this case the clopen subsets of $\Vectk^n$ are precisely the Boolean combinations of
$k$-affine subspaces.


\subsubsection{Random graph.} A random graph can be viewed as a binary relation 
$E\subset V\times V$ and thereby a f.o.structure. Its f.o.theory is an example 
of a {\em simple} theory, which is a class of f.o.theories larger and 
less well behaved then stable f.o.theories.

\subsection{A test question: can a pure Abelian group interpret a field ?} 
It is known in model theory that the first order theory of a pure Abelian group 
cannot interpret a field. Can we formulate and prove simplicially a statement 
that would imply this result in model theory ? 
Simplicially, this should mean that there are no ``non-trivial'' (in some precise sense)
surjective morphism between the corresponding f.o.theories. 
Below we define precisely the f.o.theories involved in this statement. 
\subsubsection{The f.o.theory of a pure group}  Let $(G,\cdot)$ be a group.
Equip $|G|_\bullet$ with the coarsest clopen pre-topology such that
the graph of multiplication $\cdot:G\times G\lra G$ is clopen. 
The profinite completion map 
$|G|_\bullet\lra G_\bullet$ gives its f.o.theory $G_\bullet$. 
A non-trivial theorem of Z.Sela shows that 
for a hyperbolic group $G$ its f.o.theory $G_\bullet$ is stable,
i.e. satisfies the equivalent conditions of Theorem~\ref{thm:stable}.

\subsubsection{The f.o.theory of a pure field}  Let $(F,+,\cdot)$ be a field.
Equip $|F|_\bullet$ with the coarsest clopen pre-topology such that
the graphs of addition and  multiplication $+,\cdot:F\times F\lra F$ are clopen. 
The profinite completion map 
$|F|_\bullet\lra F_\bullet$ gives its f.o.theory $F_\bullet$. 
The f.o.theory of a pure field is $\omega$-stable iff the field is  algebraically closed.
A non-trivial theorem of A.Macyntire shows
that a $\omega$-stable or superstable infinite field  is necessarily algebraically closed
\cite[\S3.1,Thm.3.1;\S6.3,Thm.6.11]{PoizatStableGroups}.

\subsubsection{A test question by Assaf Hasson} Can our simplicial methods prove that
there are no ``non-trivial'' surjective morphism $G_\bullet\lra F_\bullet$
from a pure f.o.group to a pure f.o.field ? 
It is known in model theory that a pure f.o.group does not interpret a pure f.o.field,
and an interpretation of one f.o.theory into another induces a morphism in $\sStclo$ 
or $\sStcl$. 

Assaf Hasson feels this might be an interesting test question for our simplicial approach,
even though the answer is known in model theory. The reasons are as follows.
In model theory it is hard to show
that one theory does not interpret another, not unlike in pre-algebraic topology 
it was hard to prove that there are no ``non-trivial'' maps between topological spaces.
 To give an interpretation of 
a f.o.theory $T'$ in a f.o.theory $T''$ is the same as to give a 
``non-trivial'' surjective morphism 
$T''_\bullet \lra T'_\bullet$ or $\Hommm {\sStcl} {n^\leqslant} {T''_\bullet} \lra T'_\bullet$
for some $n>0$. 


\section{A stable theory, simplicially}
Our reformulations 
are sufficient to state in a purely category-theoretic language 
a deep result in model theory. It states different characterisations of 
a well-studied notion of a {\em stable} 
f.o.theory. 

\subsection{Motivation} 
One motivation is that being stable is a necessary condition to being
uncountably categorical. Let us explain. 
A f.o.theory $T_\bullet$ is {\em categorical in cardinality $\kappa$}
iff all its models of cardinality $\kappa$ are isomorphic, 
i.e. for any two sets $M$ and $M'$ of cardinality $\kappa$ 
and any two models $M_\bullet \lra T_\bullet$ and $M'_\bullet\lra T_\bullet$ 
in $\sStclo$ there is an isomorphism $M_\bullet \xra{\approx} M'_\bullet$ over $T_\bullet$.
Morley theorem says that a f.o.theory categorical in one large enough cardinality 
is categorical in each large enough cardinality, where ``large enough''
means at least the cardinality of the language 
of the f.o.theory (= the set of all clopen subsets
of $T_\bullet(n^\leqslant)$, $n\geqslant 1$).
It is also known that such a f.o.theory is stable \cite[Ch.3]{TZ}. 

The f.o.theory of an algebraically closed field of a given characteristic 
is a representative example:
all algebraically closed fields of the same characteristic and 
of the same uncountable cardinality are isomorphic. 

\subsection{Characterisations of a f.o.theory being stable}
The following theorem states various equivalent characterisation of 
a stable first order theory. The first items $(2-4)$ are stated in purely simplicial language;
we have no nice reformulation of (5-6). 
Expressions in quotes use standard model theoretic terminology and 
are there to aid intuition but formally are unnecessary.

\begin{theo}\label{thm:stable} 
For each 
a functor $T_\bullet:\FiniteNonEmptySets\lra \Stpo$ 
satisfying eq. \eqref{TarskiVaughtRight} or eq. \eqref{TarskiVaughtLeftBC} 
the following are equivalent: 
\bee
\item ``the first order theory $T$ is stable''
\item\label{thm_stable_indi_set} ``each indiscernible sequence of $n$-tuples of a monster model of $T_\bullet$ 
is set indiscernible'', i.e. 
the following diagram stating a lifting property holds:
\bi\item $\xymatrix@C=2.39in{ 
n^\leqslant_\bullet\times \BbbQ^\leqslant_\bullet/\Aut(\BbbQ^\leqslant) \ar[r]|\forall \ar[d] & T_\bullet  \\
n^\leqslant_\bullet\times |\BbbQ|_\bullet/\SymQ \ar@{-->}[ru]|\exists & } $
\item[]
or, the same diagram in another notation, 
\item $\xymatrix@C=2.39in{ 
n^\leqslant_\bullet\times (\circlearrowleft^{e=e\circ e})_\bullet \ar[r]|\forall \ar[d] & T_\bullet \\
n^\leqslant_\bullet\times |\BbbQ|_\bullet/\SymQ \ar@{-->}[ru]|\exists & } 
$\ei
\item \label{thm_stable_few_types} ``there are few types'', 
i.e. there is a cardinal $\kappa$ (equiv., there are arbitrarily large cardinals $\kappa$)
such that for each arrow $\kappa^\leqslant_\bullet\lra T_\bullet$
there are at most $\kappa$ 
liftings 
making the diagram below commute
\\ \xymatrix@C=2.39in{ 
& T_\bullet\circ[+1] \ar[d] \\
\kappa^\leqslant_\bullet \ar[r] \ar@{-->}[ru]  & T_\bullet } 
\item ``types over models are definable'', i.e.
each type over an $<\!\!\aleph_0$-saturated model $N$ is definable over $N$, 
i.e. for the theory $T(N)$ denoting the theory $T$ extended by constants $N$, it holds
(where the lifting is required to be continuous)
\\ \xymatrix@C=2.29in{ |N|_\bullet \ar[r] \ar[d]|{(N\text{ is a model }T_\bullet)} 
& T(N)_\bullet\circ[+1] \ar[d] \\
T(N)_\bullet \ar[r]|\id \ar@{-->}[ru]  & T(N)_\bullet } 
\item \label{stable:orderpr} 
``no formula $\varphi(x,y)$ has the order property'', i.e. 
there is no formula $\varphi(x,y)$ such that for each $N$
there is a sequence $(a_i,b_i)_{0<i< N}$ of tuples in some, equiv. any, model $M$ 
\bi\item[] $M\models\varphi(a_i,a_j) \,\,\Longleftrightarrow\,\, i\leqslant j$\ei
\item \label{stable:Gro} in an $\leqslant\aleph_0$-saturated model for each formula $\varphi(x,y)$ 
and each sequence $(a_i,b_i)_{i\in\Bbb N}$ of tuples it holds
$$\lim\limits_{i\rightarrow\infty}\lim\limits_{j\rightarrow\infty}\varphi(a_i,b_j)
=
\lim\limits_{j\rightarrow\infty}\lim\limits_{i\rightarrow\infty}\varphi(a_i,b_j)$$
\eee
\end{theo}
\begin{proof} (1)$\iff$(2)$\iff$(3)$\iff$(4)$\iff$(5):
This is a standard theorem in model theory, 
e.g. see the textbook \cite[Thm 8.2.3, Exercise 8.2.6-7,Thm 8.3.1]{TZ}. 
\newline\ref{stable:orderpr}$\iff $\ref{stable:Gro} is evident.
The latter form of the order property is an instance of a compactness criterion defined by  
Grothendieck in analysis in terms of commuting limits 
\cite{GrothendieckStable,benYaacov}. %
\end{proof} 
\begin{rema}
For definability of types 
I don't have a nice description of $T(N)_\bullet$.
\end{rema}

\begin{defi}[Stable f.o.theory] A f.o.theory $T_\bullet$ is {\em stable} iff it satisfies
the equivalent conditions of Theorem~\ref{thm:stable}. 
A f.o.theory is {\em $\kappa$-stable} iff it satisfies Theorem~\ref{thm:stable}\ref{thm_stable_few_types}(few types) for the cardinal $\kappa$. 
A f.o.theory is {\em superstable} iff it is $\kappa$-stable for all $\kappa$ large enough.
\end{defi}

\subsection{A stable f.o.theory defined in term of forking/independence relation}
Algebraic independence is an example of an independence relation in the case 
of algebraically closed field, and so is linear independence in vector fields. 
An independence relation with similar properties exists for any stable theory,
and the notion of a stable theory can be characterised in terms of existence 
of such an independence relation, which, moreover, is necessarily unique. 
\subsubsection{A model theoretic informal explanation}
For subsets $A,B,C\subset \Bbb C$ of $\Bbb C$,
we may say that $A$ is algebraically, resp.~$\Bbb Q$-linearly,  independent from $B$ over $C$
and write  $  A \forkindep[C] B $ iff for each finite tuple $a_1,..,a_n\in A$ 
the transendence degree, resp.~linear dimension over $\Bbb Q$, 
of $a_1,...,a_n$ over $B$ is equal to that over $C$.  
In model theory for any f.o.theory $T$ Shelah defines 
a {\em non-forking} relation $  A \forkindep[C] B $ on subsets $A,B,C$ of a monster model $M$.
Stable theories can be characterised as those where this relation is particularly well-behaved
in a way reminiscent of the algebraic or linear independence relation. 
It is convinient to view an independence relation $  a \forkindep[C] B $ on subsets 
$a,C,B\subset M$ of a model 
as  a notion (=distinguished class) of a free  extension of types which in model theory 
is known as a {\em non-forking extension of types}: for a finite tuple $a$ and sets $B,C$ 
we say the type $\tp(a/BC)$ of tuple $a$ with parameters in $B\cup C$ is 
an {\em non-forking} extension of the type $\tp(a/C)$ iff $a$ is independent from $B$ over $C$. 
We write
\bi\item $   a \forkindep[C] B $ iff $\tp(a/C) \sqsubset \tp(a/BC)$.\ei
Recall that {\em the type $\tp(a/C)$ of tuple $a$ over a set $C$ of parameters} is defined
to be the set of all first order  formulas with parameters in $C$ true of $a$. 
The same data are represented by the commutative triangle (where $a=(a_1,...,a_n)$) 
$$\xymatrix@C=2.29in{ & T\circ[+n] \ar[d] \\ 
|C|_\bullet \ar[r]^{(c_1,...,c_l)\mapsto \tp(c_1,...,c_l)} \ar[ru]|{\tp(a/C)} & T_\bullet }$$

\subsubsection{A characterisation of stability in terms 
of a distinguished class of lifting diagrams}
The theorem quotes the axioms of a binary relation $p\sqsubset q$ on types as stated 
in \cite[p.246]{HarnikHarrington}, and then expresses them simplicially. 
Thus, the reader does not need to understand the expressions in quotes.
A more direct simplicial translation should 
replace everywhere below  $\alpha^\leqslant_\bullet$ by $|\alpha|_\bullet$ 
and work in the category of symmetric simplicial sets. 

Note that the commutative diagrams in the theorem may be considered 
in the category of simplicial sets as the topology is not used. 
Perhaps it can be speculated that this class of distinguished diagrams plays 
the role of topology. 

\begin{theo}\label{thm:stableViaForking}
A f.o.theory is stable iff there is a collection of distinguished commutative triangles
of the form below where the arrow $\alpha^\leqslant_\bullet\lra\beta^\leqslant_\bullet$ is an arbitrary injective morphism between simplicial sets represented by cardinals $\alpha$ and $\beta$
$$\xymatrix@C=2.39in{
\alpha^\leqslant_\bullet \ar[r]|p \ar@{>->}[d] & T_\bullet\circ[+n] 
\ar[d]|{\pr^{n+1,n+2,...}_\bullet}\\
\beta^\leqslant_\bullet \ar[r] \ar@{->}[ru]|q ^{p\sqsubset q}  & T_\bullet
}$$
satisfying the following Axioms:
\bee
\item[Axiom 2.] 
``If $p \in S(A)$ and $A \subset  B$, then $p\sqsubset q$ for some $q \subset S(B)$.''
$$\xymatrix@C=2.39in{
\alpha^\leqslant_\bullet \ar[r]|p \ar@{>->}[d] & T_\bullet\circ[+n] 
\ar[d]|{\pr^{n+1,n+2,...}_\bullet}\\
\beta^\leqslant_\bullet \ar[r] \ar@{-->}[ru]|{\exists q} ^{(p\sqsubset q)}  & T_\bullet
}$$

\item[Axiom 1.] 
``If $p \sqsubset q \sqsubset r$, then:
(a) $p\sqsubset q \sqsubset r$ implies that $p \sqsubset  r$''
$$
\xymatrix@C=2.39in{
\alpha^\leqslant_\bullet \ar[r]|p \ar@{>->}[d] & T_\bullet\circ[+n] 
\ar[dd]|{\pr^{n+1,n+2,...}_\bullet}\\
\beta^\leqslant_\bullet \ar@{>->}[d] 
\ar@{->}[ru]|q^{(p\sqsubset q)}_{(q\sqsubset r)} & \\
\gamma^\leqslant_\bullet \ar[r] \ar@{->}[ruu]|r_{\therefore(p\sqsubset r)} & T_\bullet
}$$
``If $p \sqsubset q \sqsubset r$, then:
(b) $p \sqsubset r$ implies that $p \sqsubset q$,
(c) $p \sqsubset r$ implies that $q\sqsubset r$.''
$$
\xymatrix@C=2.39in{
\alpha^\leqslant_\bullet \ar[r]|p \ar@{>->}[d] & T_\bullet\circ[+n] 
\ar[dd]|{\pr^{n+1,n+2,...}_\bullet}\\
\beta^\leqslant_\bullet \ar@{>->}[d] 
\ar@{->}[ru]|q^{\therefore(p\sqsubset q)}_{\therefore(q\sqsubset r)} & \\
\gamma^\leqslant_\bullet \ar[r] \ar@{->}[ruu]|r_{(p\sqsubset r)} & T_\bullet
}$$
\item[Axiom 4.] ``For any $p$ there is a cardinal $\lambda$ such that there are 
at most $\lambda$ mutually contradictory types $q$ s.t. $p \sqsubset q$.''
\\ For any $p:\alpha^\leqslant_\bullet\lra T_\bullet[+n]$ 
there is a cardinal $\lambda$ such that 
for any cardinal $\beta>\alpha$ 
for each commutative square as shown below 
there are 
at most $\lambda$ 
liftings $q:\beta^\leqslant_\bullet \lra T_\bullet[+n]$ s.t. $p \sqsubset q$.
$$
\xymatrix@C=2.39in{
\alpha^\leqslant_\bullet \ar[r]|p \ar@{>->}[d] & T_\bullet\circ[+n] 
\ar[d]|{\pr^{n+1,n+2,...}_\bullet}\\
\beta^\leqslant_\bullet \ar[r] \ar@{-->}[ru]|{\exists^{\leqslant\lambda} q \,\,(p\sqsubset q)}  & T_\bullet
}$$
or, equivalently, \\ for each $p:\alpha^\leqslant_\bullet\lra T_\bullet\circ[+n]$
there is  a cardinal $\lambda$ such that among any $>\lambda$ 
morphisms $q:\beta^\leqslant_\bullet \lra T_\bullet\circ[+n]$ such that $p\sqsubset q$
there are at least two morphisms 
$q':\beta'^\leqslant_\bullet \lra T_\bullet\circ[+n]$ and $q'':\beta''^\leqslant_\bullet \lra T_\bullet\circ[+n]$ admitting a factorisation though 
some $r:\gamma^\leqslant_\bullet \lra T_\bullet\circ[+n]$
$$
\xymatrix@C=3.39cm{
&\alpha^\leqslant_\bullet \ar[rr]|p \ar@{>->}[dl] \ar@{>->}[dr]  & & T_\bullet\circ[+n] 
\ar[d]^{\pr^{n+1,n+2,...}_\bullet}\\
\beta'^\leqslant_\bullet  \ar[rrru]|{(p\sqsubset q')} \ar@{>-->}[dr]|{\exists} 
\ar@/^1pc/[rrr] &&
\beta''^\leqslant_\bullet  \ar[ru]|{(p\sqsubset q'')} \ar@{>-->}[dl]|{\exists} \ar[r] & T_\bullet \\
& \gamma^\leqslant_\bullet \ar@{-->}@/_2pc/[rru]|\exists  \ar@{-->}@/_2pc/[rruu]|\exists
}$$

\item[Axiom 3.] 
``There is a cardinal $\kappa$ such that if $p \in S(A)$, then $p | A_0 \sqsubset p$ 
for some $A_0 \subset  A$ with $\card (A_0) < \kappa$.''
\\ There is a cardinal $\kappa$ such that
for each cardinal $\beta$ and each $p:\beta^\leqslant_\bullet\lra T_\bullet[+n]$
there is cardinal $\alpha<\kappa$, an injective monotone map 
$\alpha\lra \beta$, and $p':\alpha ^\leqslant_\bullet\lra T_\bullet[+n]$
such that $p'\sqsubset p$ making the diagram below commute.
$$\xymatrix@C=2.39in{
\alpha^\leqslant_\bullet \ar@{-->}[r]|{\exists p'} 
\ar@{>-->}[d]_{\exists}^{\alpha<\kappa } & T_\bullet\circ[+n] 
\ar[d]|{\pr^{n+1,n+2,...}_\bullet}\\
\beta^\leqslant_\bullet \ar[r] \ar@{->}[ru]|{p}^{(p'\sqsubset p)}   & T_\bullet
}$$
\eee
\end{theo} 
\begin{proof} This is a standard result, 
see e.g.~\cite[Introduction;Theorem 5.8]{HarnikHarrington}.\end{proof} 

\begin{rema} A f.o.theory is superstable iff it satisfies the Axioms above 
and the Axiom 3 with $\kappa:=\aleph_0$.
\bee
\item[Axiom $3_{\aleph_0}$.] ``If $p \in S(A)$, then $p | A_0\sqsubset p$ 
for some finite $A_0\subset A$.''
 This is the Axiom 3 for $\kappa:=\aleph_0$. `?todo: It is tempting to modify our definitions
so that it can expressed using $T_\bullet\circ[+m]\lra T_\bullet$.?
$$\xymatrix@C=2.39in{
m^\leqslant_\bullet \ar@{-->}[r]|{\exists p'}_{(p'\sqsubset p)}
\ar@{>-->}[d]|{\exists \alpha<\kappa } & T_\bullet\circ[+n] 
\ar[d]|{\pr^{n+1,n+2,...}_\bullet}\\
\beta^\leqslant_\bullet \ar[r] \ar@{->}[ru]|{p}   & T_\bullet
}$$
\eee
\end{rema} 

\subsubsection{Non-dividing extension of types} There is an explicit way 
to define a freeness property of extension of types which under some assumptions
is equivalent to non-dividing, but in general stronger and less well behaved.
It uses an object used in Theorem~\ref{thm:stable}\ref{thm_stable_indi_set}. 
We formulate the definition only for types over the empty set. 
It is a reformulation of \cite[Cor.7.1.5(2)]{TZ}.
\begin{defi} A type $tp(a_1,..,a_m/b_1,..,b_n)$ {\em does not divide over $\emptyset$} iff
\bi
\item $\xymatrix@C=2.39in{
n^\leqslant_\bullet \ar[d] \ar[r]|{\tp(a_1,...,a_m/b_1,...,b_n)} & T_\bullet[+m] 
\ar[d]|{\pr_\bullet^{m+1,m+2,...}} \\
n^\leqslant_\bullet\times \BbbQ^\leqslant_\bullet/\Aut(\BbbQ^\leqslant)
\ar[r]|\forall  \ar@{-->}[ru]|\exists
& T_\bullet}$
\item[] or, the same diagram in another notation, 
\item $\xymatrix@C=2.39in{
n^\leqslant_\bullet \ar[d] \ar[r]|{\tp(a_1,...,a_m/b_1,...,b_n)} & T_\bullet[+m] 
\ar[d]|{\pr_\bullet^{m+1,m+2,...}} \\
n^\leqslant_\bullet\times (\circlearrowleft^{e=e\circ e})_\bullet 
\ar[r]|\forall  \ar@{-->}[ru]|\exists
& T_\bullet}$
\ei
\end{defi}

%
%


\subsection{Questions about stability} Could this reformulation of a standard basic result be of use ? We state a couple of questions it suggests. 

\begin{que}[Non-symmetric stable f.o.theories]  Does the standard model theoretic methods prove 
in fact proves the equivalence of various characterisations of stability
for  an arbitrary non-symmetric (i.e. usual) 
simplicial space $T_\bullet:\Dop\lra \Stpo$ 
satisfying \eqref{TarskiVaughtRight} or \eqref{TarskiVaughtLeftBC} ?
\end{que} 

\begin{que}[A category theoretic proof and reformulation]
 Do characterisations of stability  remind of standard results 
 in category theory ? What are precise assumptions on the category 
 $\Stpo$ and the functor $T_\bullet:\Dop\lra \Stpo$ used to prove the equivalence 
 of various characterisations ? 
\end{que}

\begin{rema} This reformulation is unsatisfactory in many obvious ways. 
Its only satisfactory aspect is purely logical rather than mathematical: 
it is very {\em tautological}: merely a change of notation expressing by 
standard simplicial diagrams standard sentences in model theory.
\end{rema}

\section{Standard category-theoretic constructions applied to a f.o.theories}
In this section we make a couple of definitions which would be hard to formulate
in purely model theoretic terminology.  This section is preliminary and relies 
on future work. 

\subsection{Standard universal constructions using the lifting property}
Category theory allows us to define a number of universal constructions.

We shall now suggest a couple of definitions based on 
a standard diagram chasing trick arising in homotopy theory
called the lifting property or weak orthogonality of morphisms%
\footnote{Recall that a morphism $i$ in 
a category has the {\em left lifting property} with respect 
to a morphism $p$, 
and $p$ also has the {\em right lifting property} with respect to $i$, 
denoted $i\rtt p$, 
iff for each $f:A\to X$ and $g:B\to Y$ such that $p\circ f = g \circ i$ there exists $h:B\
to X$ such that $h\circ i = f$ and $p\circ h = g$.

For a class $P$ of morphisms in a category, its {\em left orthogonal} $P^{\rtt l}$ with respect to the lifting property, respectively its {\em right orthogonal} $P^{\rtt r}$, is the class of all morphisms which have the left, respectively right, lifting property with respect to each morphism in the class $P$. In notation,
$$
P^{\rtt l} := \{ i \,\,:\,\, \forall p\in P\,\, i\rtt p\},
P^{\rtt r} := \{ p \,\,:\,\, \forall i\in P\,\, i\rtt p\}, P^\lr:=(P^\lrl)^\rlr,..$$

Taking the orthogonal of a class $P$ is a simple way to define a class of morphisms excluding non-isomorphisms 
from $P$, in a way which is useful in a diagram chasing computation, and is often used to
define properties of morphisms starting from an explicitly given class of (counter)examples. 
A number of standard basic notions may also be expressed using the lifting property starting from a list of (counter)examples such as a group being nilpotent, soluble, free, abelian, $p$-group, or a topological space being contractible, compact, having a generic point, and others. 
  For this reason, 
  it is convenient and intuitive to refer to $P^\lrl$ and $P^\rlr$ 
  as {\em left, resp.~right, Quillen negation} of property $P$.}
used in our reformulations.

\subsection{A notion of $T$-like f.o.theory} 
For any property (=class) $P$ of morphisms we may define its {\em $lr$-generalisation 
$P^{lr}\supset P$}. Pick a f.o.theory $T_\bullet$ and take 
$P_T:=\{ T_\bullet\lra\top\}$. 
The notion of a f.o.theory is defined by a right lifting property 
\eqref{TarskiVaughtRight}, and therefore 
$Q_\bullet\lra\top \,\in\, \big\{\, T_\bullet\lra\top\,\big\}^{lr}$ implies 
that $Q_\bullet$ is also a f.o.theory for a f.o.theory $T_\bullet$.

It may be useful to enlarge the class (= the property of being $T_\bullet$) 
$\{ T_\bullet\lra \top\}$ by adding some morphisms satisfying \eqref{TarskiVaughtRight},
e.g. the morphism from the definition \eqref{TarskiVaughtLeft} of a model of a f.o.theory.  That is, to consider 
$$P'_T:=\big\{\, T_\bullet \lra \top\,,\,\, 
|\{\star_1,\star_2\}|_\bullet \sqcup |\{\star_1\}|_\bullet\lra
|\{\star_1,\star_2\}|_\bullet \,  \big\}$$

This leads us to the following definitions.

\begin{defi}[$T_\bullet$-like] Say that a f.o.theory $Q_\bullet$ is {\em $T_\bullet$-like} iff 
\begin{equation} {{Q_\bullet}\atop{\downarrow\atop\top}} \in
\left\{ {{T_\bullet}\atop{\downarrow\atop\top}} \,\,\,\,
{ {{ |\{\star_1,\star_2\}|_\bullet \sqcup |\{\star_1\}|_\bullet }}
\atop{\downarrow\atop
{ \displaystyle |\{\star_1,\star_2\}|_\bullet \phantom{\sqcup |\{\star_1\}|_\bullet}  }}}
\right\}^{lr}\end{equation}
or, dually, 
\begin{equation} {{\bot}\atop{\downarrow\atop\displaystyle Q_\bullet}} \in
\left\{ {{\bot}\atop{\downarrow\atop\displaystyle T_\bullet}} \,\,\,\,
{ { (k+m)^\leqslant_\bullet\vee_{k^\leqslant_\bullet} (k+n)^\leqslant_\bullet }
\atop{\downarrow\atop
{\displaystyle  (k+m+n)^\leqslant_\bullet }}}
\right\}^{rl}\end{equation}
\end{defi}

\begin{defi}[$P$-interpolant] For a class (=property) of morhisms, say 
that a f.o.theory $Q_\bullet$ is a {\em $P$-$l.lr$-interpolant of f.o.theory $T'_\bullet$
and $T''_\bullet$} iff there is a decomposition
\begin{equation}T'_\bullet \xra {(P)^l} Q_\bullet \xra{(P)^{lr}} T''_\bullet
\end{equation}
\end{defi}

\begin{que} Do the notions of $T$-like and $P$-$l.lr$-interpolant 
have model theoretic meaning ? 
\end{que} 
\begin{proof}[Considerations] Say, it seems that  $T$-like f.o.theories
inherit some of the properties of generalised indiscernibles, as those 
are sometimes expressed by left lifting properties. 
\end{proof} 

\begin{que} What can we say about the properties of theories defined by $P_T^{lr}$
or ${P'}_T^{lr}$ ? 
\end{que} 

\subsection{Decalage cohomology of a f.o.theory ?} This subsection relies on a notion defined in detail in future work. 
Can we use 
decalage cohomology 
to define non-trivial 
invariants of f.o.theories ? 
Note that our reformulation has as instances the standard 
chain complexes of Cech cohomology of a sheaf wrt an open covering 
and of group cohomology 
and thus likely has the same functorial properties, 
e.g. long exact sequence of cohomology associated with something like 
a short exact sequence.

\subsubsection{A definition of cohomology via decalage} 
Let $\PPhi$ be a category, and 
let $\sPPhi:=Func(\Dop,\PPhi)$ be its category of simplicial objects. 

Let $p_\bullet:A_\bullet \lra B_\bullet$ be an abelian group object 
in the slice category $\sPPhi_{/B_\bullet}$ over $B_\bullet:\sPPhi$. 

Use decalage (shift)
$[+1]:\Dop\lra\Dop, n^\leqslant\longmapsto (1+n)^\leqslant$ 
to define a cosimplicial set 
as follows:
\begin{equation}\label{CoSimpDec}
\calA^{n-1}\left({{A_\bullet}\atop{\downarrow p_\bullet\atop{B_\bullet}}}\right) 
:= \Hommm {\sPPhi_{/B_\bullet}} { B_\bullet\circ[+n] 
\xra {\pr_{1,..,n,\bullet}} B_\bullet\text{  }} 
{\text{  }A_\bullet\xra{p_\bullet} B_\bullet}
\end{equation}

In other words, define {\em an $n-1$-cochain} to be a commutative triangle in $\sPPhi$
\begin{equation}\label{CoSimpDecTriangle}
\xymatrix@C=2.39cm{ & A_\bullet \ar[d] \\
B_\bullet\circ[+n] \ar[r]|{\pr_{1,..,n,\bullet}}  \ar[ru]|{\sigma_\bullet} & B_\bullet }
\end{equation}
By funtoriality $A_\bullet\lra B_\bullet$ being an abelian group object in 
$\sPPhi_{/B_\bullet}$ 
implies that $\calA^\bullet$ defines a cosimplicial abelian group.
That is, a map $+:A_\bullet\times_{B_\bullet} A_\bullet \lra A_\bullet$ over $B_\bullet$
determines in an obvious way for each $n\geqslant 0$ a map $$
+: \calA^{n}\left({{A_\bullet}\atop{\downarrow p_\bullet\atop{B_\bullet}}}\right) 
\times
\calA^{n}\left({{A_\bullet}\atop{\downarrow p_\bullet\atop{B_\bullet}}}\right) 
\lra \calA^{n}\left({{A_\bullet}\atop{\downarrow p_\bullet\atop{B_\bullet}}}\right) 
$$
\begin{equation}\label{CoSimpDecTrianglePlus}
\xymatrix@C=2.39cm{ & A_\bullet \ar[d] \\
B_\bullet\circ[+n+1] \ar[r]|>>{\pr_{1,..,n+1,\bullet}}  
\ar[ru]|{\sigma'_\bullet+\sigma''_\bullet} & B_\bullet }
\end{equation}

By the standard Dold-Kan
correspondence this gives a chain complex and therefore a cohomology theory we denote by
$H^\bullet(B_\bullet, {{A_\bullet}\atop{\downarrow p_\bullet\atop{B_\bullet}}})$. 

\begin{defi}[Decalage cohomology]\label{defi:DefCoh} For a cosimplicial abelian group 
object 
$p_\bullet:A_\bullet\lra B_\bullet$ over $B_\bullet$ in a category $\sPPhi_{/B_\bullet}$, 
 {\em the decalage cohomology 
 $H^\bullet(B_\bullet,\,{{A_\bullet}\atop{\downarrow p_\bullet\atop{B_\bullet}}})$ 
 of $B_\bullet$ with coefficients $p:A_\bullet\lra B_\bullet$}
 is the cohomology associated with the 
the cosimplicial abelian group defined in \eqref{CoSimpDecTrianglePlus}.
\end{defi}

\subsubsection{Cohomology of a f.o.theory with constant coefficients} 
Note that Definition ~\ref{defi:DefCoh} 
makes well-defined 
the notion of 
``the decalage cohomology of a first order theory with constant integer coefficients''.
\begin{que} For a f.o.theory $T_\bullet$, is 
its decalage cohomology with constant integer coefficients necessarily trivial ? 
$$H(T_\bullet, \constb \Bbb Z\times T_\bullet \lra T_\bullet)=0 ? $$ 
Say, can decalage cohomology (with what coefficients?) distinguish
the stable f.o.theory
$\faktour{\BbbC_\bullet}{\AutCQ}$ from the unstable profinite compactification of 
$\faktour{\BbbR_\bullet}{\AutRQ}$ ? 
\end{que}

\subsubsection{Cohomology of a type or a definable type} 
We could modify our definition of decalage cohomology to 
define cohomology of a type, as follows.

Recall that a $l$-type $p(v_1,...,v_l/A)$ of theory $T$ over a set $A$
is the same as a map $p_\bullet:|A|_\bullet\lra T_\bullet\circ[+l]$, 
and that a $\emptyset$-definable $n$-type is the same as a map 
$p_\bullet:T_\bullet\lra T_\bullet\circ[+l]$. 

\begin{que} 
Does the following definition define a non-trivial model theoretic property of a type,
a definable type, or an indiscernible sequence ? 

For a map $p_\bullet:C_\bullet \lra T_\bullet\circ[+1]$ define its decalage cohomology with constant coefficients $\Bbb Z$ by defining a $n$-cochain to be a commutative triangle
\begin{equation}\label{CoSimpDecTriangle:p_type}
\xymatrix@C=2.39cm{ && \Bbb Z \times T_\bullet\circ[+1] \ar[d] \\
C_\bullet\circ[+n] \ar[r]|{\pr_{1,..,n-1,\bullet}} \ar[r]  
\ar[rru]|{\sigma_\bullet} & C_\bullet \ar[r]|{p_\bullet} & T_\bullet\circ[+1] }
\end{equation} 

Cohomology of a definable type $p_\bullet:T_\bullet\lra T_\bullet\circ[+1]$ 
would be defined by cochains
\begin{equation}\label{CoSimpDecTriangle:p_deftype}
\xymatrix@C=2.39cm{ && \Bbb Z \times T_\bullet\circ[+1] \ar[d] \\
T_\bullet \circ[+n] \ar[r]|{\pr_{1,..,n-1,\bullet}} \ar[r]  
\ar[rru]|{\sigma_\bullet} & T_\bullet  \ar[r]|{p_\bullet} & T_\bullet\circ[+1] }
\end{equation} 
\end{que}
Cohomology of a (usual) type $p_\bullet:|A|_\bullet\lra T_\bullet\circ[+1]$ over a set $A$ of parameters would be defined by cochains
\begin{equation}\label{CoSimpDecTriangle:p_deftype}
\xymatrix@C=2.39cm{ && \Bbb Z \times T_\bullet\circ[+1] \ar[d] \\
|A|_\bullet \circ[+n] \ar[r]|{\pr_{1,..,n-1,\bullet}} \ar[r]  
\ar[rru]|{\sigma_\bullet} & |A|_\bullet  \ar[r]|{p_\bullet} & T_\bullet\circ[+1] }
\end{equation}

For a complete EM-type of an indiscernible sequence 
$p_\bullet:(\circlearrowleft^{e=e\circ e})_\bullet \lra T_\bullet$
the cochains would be
\begin{equation}\label{CoSimpDecTriangle:p_EMtype}
\xymatrix@C=2.39cm{ && \Bbb Z \times T_\bullet \ar[d] \\
(\circlearrowleft^{e=e\circ e})_\bullet \circ[+n] \ar[r]|{\pr_{1,..,n-1,\bullet}} \ar[r]  
\ar[rru]|{\sigma_\bullet} & (\circlearrowleft^{e=e\circ e})_\bullet  \ar[r]|{\text{indi.seq. }p_\bullet} & T_\bullet }
\end{equation}

\end{document}